\RequirePackage[leqno]{amsmath}
\documentclass[leqno]{amsart}
\usepackage{color}
\usepackage{amsfonts}
\usepackage{enumitem}
\setlist[itemize]{topsep=0pt,after=\vspace{1.5\baselineskip}}
 \usepackage[colorlinks=true]{hyperref}
\usepackage{empheq}
\usepackage[T1]{fontenc}

\def\R{\mathbb R}  
\def
\@cite
#1#2{[{{\bfseries #1}\if@tempswa , #2\fi}]}
\newtheorem{theorem}{Theorem}[section]

\newtheorem{lemma}[theorem]{Lemma}

\newtheorem{remark}{Remark}

\title[{\small{B\MakeLowercase{low-up prevention in an attraction--repulsion
chemotaxis system}}}] 
     {\large{B\MakeLowercase{low-up prevention by sublinear production in a $n$-dimensional attraction-repulsion chemotaxis system}}}
\author[N. Pintus and G. Viglialoro]{}
\keywords{Chemotaxis, boundedness, blow-up prevention.\\
\textit{$^\star$Corresponding author}: giuseppe.viglialoro@unica.it\\
\textit{Acknowledgements.} The authors are members of the Gruppo Nazionale per l'Analisi Matematica, la Probabilit\`a e le loro Applicazioni (GNAMPA) of the Istituto Na\-zio\-na\-le di Alta Matematica (INdAM). GV is partially supported by the research projects \textit{Integro-differential Equations and Non-Local Problems}, funded by Fondazione di Sardegna (2017) and MIUR (Italian Ministry of Education, University and Research) Prin 2017 \textit{Nonlinear Differential Problems via Variational, Topological and Set-valued Methods} (Grant Number: 2017AYM8XW)}
\RequirePackage[normalem]{ulem} 
\RequirePackage{color}\definecolor{RED}{rgb}{1,0,0}\definecolor{BLUE}{rgb}{0,0,1} 
\providecommand{\DIFaddbegin}{} 
\providecommand{\DIFaddend}{} 
\providecommand{\DIFdelbegin}{} 
\providecommand{\DIFdelend}{} 
\providecommand{\DIFaddbeginFL}{} 
\providecommand{\DIFaddendFL}{} 
\providecommand{\DIFdelbeginFL}{} 
\providecommand{\DIFdelendFL}{} 
\newcommand{\DIFscaledelfig}{0.5}
\RequirePackage{settobox} 
\RequirePackage{letltxmacro} 
\newsavebox{\DIFdelgraphicsbox} 
\newlength{\DIFdelgraphicswidth} 
\newlength{\DIFdelgraphicsheight} 
\LetLtxMacro{\DIFOincludegraphics}{\includegraphics} 
\newcommand{\DIFaddincludegraphics}[2][]{{\color{blue}\fbox{\DIFOincludegraphics[#1]{#2}}}} 
\newcommand{\DIFdelincludegraphics}[2][]{
\sbox{\DIFdelgraphicsbox}{\DIFOincludegraphics[#1]{#2}}
\settoboxwidth{\DIFdelgraphicswidth}{\DIFdelgraphicsbox} 
\settoboxtotalheight{\DIFdelgraphicsheight}{\DIFdelgraphicsbox} 
\scalebox{\DIFscaledelfig}{
\parbox[b]{\DIFdelgraphicswidth}{\usebox{\DIFdelgraphicsbox}\\[-\baselineskip] \rule{\DIFdelgraphicswidth}{0em}}\llap{\resizebox{\DIFdelgraphicswidth}{\DIFdelgraphicsheight}{
\setlength{\unitlength}{\DIFdelgraphicswidth}
\begin{picture}(1,1)
\thicklines\linethickness{2pt} 
{\color[rgb]{1,0,0}\put(0,0){\framebox(1,1){}}}
{\color[rgb]{1,0,0}\put(0,0){\line( 1,1){1}}}
{\color[rgb]{1,0,0}\put(0,1){\line(1,-1){1}}}
\end{picture}
}\hspace*{3pt}}} 
} 
\LetLtxMacro{\DIFOaddbegin}{\DIFaddbegin} 
\LetLtxMacro{\DIFOaddend}{\DIFaddend} 
\LetLtxMacro{\DIFOdelbegin}{\DIFdelbegin} 
\LetLtxMacro{\DIFOdelend}{\DIFdelend} 
\DeclareRobustCommand{\DIFaddbegin}{\DIFOaddbegin \let\includegraphics\DIFaddincludegraphics} 
\DeclareRobustCommand{\DIFaddend}{\DIFOaddend \let\includegraphics\DIFOincludegraphics} 
\DeclareRobustCommand{\DIFdelbegin}{\DIFOdelbegin \let\includegraphics\DIFdelincludegraphics} 
\DeclareRobustCommand{\DIFdelend}{\DIFOaddend \let\includegraphics\DIFOincludegraphics} 
\LetLtxMacro{\DIFOaddbeginFL}{\DIFaddbeginFL} 
\LetLtxMacro{\DIFOaddendFL}{\DIFaddendFL} 
\LetLtxMacro{\DIFOdelbeginFL}{\DIFdelbeginFL} 
\LetLtxMacro{\DIFOdelendFL}{\DIFdelendFL} 
\DeclareRobustCommand{\DIFaddbeginFL}{\DIFOaddbeginFL \let\includegraphics\DIFaddincludegraphics} 
\DeclareRobustCommand{\DIFaddendFL}{\DIFOaddendFL \let\includegraphics\DIFOincludegraphics} 
\DeclareRobustCommand{\DIFdelbeginFL}{\DIFOdelbeginFL \let\includegraphics\DIFdelincludegraphics} 
\DeclareRobustCommand{\DIFdelendFL}{\DIFOaddendFL \let\includegraphics\DIFOincludegraphics} 

\begin{document}
\maketitle
\DIFdelbegin 
\DIFdelend 

\centerline{\scshape Nicola Pintus$^\sharp$ and Giuseppe  Viglialoro$^{\sharp,\star}$}
\medskip
{\footnotesize
 \centerline{$^\sharp$Dipartimento di Matematica e Informatica}
 \centerline{Universit\`{a} di Cagliari}
 \centerline{V. le Merello 92, 09123. Cagliari (Italy)}
}

\bigskip
\begin{abstract}
In this paper we study a zero-flux attraction-repulsion che\-mo\-taxis-system. We show that despite any mutual interplay between the repulsive and attractive coefficients from the corresponding chemo-sensitivities, even less any restriction on their own sizes, if the production rate of that chemical signal responsible of the cellular coalescence is sublinear, then any initial data emanate a unique global  classical solution, which is as well bounded.  Further, in a remark of the manuscript, we also address an open question given in \cite{Viglialoro2019RepulsionAttraction}.

\end{abstract}
\section{Introduction and motivations: presentation of the main theorem}\label{IntroductionSection} 
The biological models presented by Keller and Segel  in their landmarking papers \cite{K-S-1970,Keller-1971-MC}, and describing chemotaxis phenomena, have  been lately inspiring many interesting investigations in the fields of both theoretical and applied mathematics.  In this sense, aim of the present research is focusing on a precise variant (which, as technically detailed in \cite{Luca2003Alzheimer}, fits with real applications concerning aggregation phenomena of microglia observed in Alzheimer's disease) of the aforementioned models  and enhance its underlying  mathematical theory. 

To be precise, this paper deals with the analysis of the mathematical problem idealizing the motion of a certain cell density $u(x, t)$ at the position $x$ and at the time $t$, initially distributed according to the law of $u_0(x):=u(x,0)$, and moving in an insulated domain under a repulsion effect, from a certain chemical signal concentration $w(x, t)$, and an attraction impact, from another one $v(x, t)$, which is ``slightly'' weaker (as specified later) than the first. Mathematically, we will face this system
\begin{equation}\label{problem}
\begin{cases}
u_{ t}=\Delta u -\chi \nabla \cdot (u\nabla v)+\xi \nabla \cdot (u\nabla w) & \textrm{in }\Omega, t>0, \\
0=\Delta v+\alpha u^\rho-\beta v & \textrm{in } \Omega,  t>0,\\
0=\Delta w+\gamma u-\delta w & \textrm{in } \Omega,  t>0,\\
 u_{\nu}= v_{\nu}= w_{\nu}=0 & \textrm{on } \partial \Omega, t>0, \\
u_0(x):=u(x,0)\geq 0 & x \in  \Omega,
\end{cases}
\end{equation}
defined in $\Omega$, a bounded and smooth domain of $\R^n$ with $n\geq 2$, and where $\alpha,\beta,$ $\gamma,\delta,\chi,\xi >0,$ $0<\rho<1$,  whilst ${(\cdot)}_{\nu}$ indicates the outward normal derivative on $\partial \Omega$ and $u_0(x)$ is a nonnegative and sufficiently regular function over $\bar{\Omega}$. 
 Physically, the parameters $\chi$ and $\xi$ measure the influences of the attraction and repulsion, whereas the second and third equations idealize that chemoattractant and chemorepellent, $v$ and $w$, are sublinearly and linearly released by cells, moreover decaying with rates $\beta$ and $\delta$, respectively.

Similarly to what happens in the abundantly studied parabolic-elliptic Keller--Segel system, obtained by \eqref{problem} letting $\xi=0$, $\rho=1$ and eliminating the third unknown $w$, also model \eqref{problem} itself is likely to manifest  the so called \textit{chemotactic collapse}, the mechanism according to which the movement of the cells may eventually degenerate into aggregation processes giving rise to $\delta$-formations. 

As far as this research is concerned, our accurate bibliographic research did not show any result about prototypes as that in \eqref{problem} presenting sublinear production. Subsequently, since we establish here that such weakening in the rate growth of the chemical signal associated to the cells' gathering suffices to prevent any impulsive concentration of the same cells independently by the coalescence effects coming from other factors of the model, we believe that this work provides a more complete picture about attraction-repulsion chemotaxis problems. 

Conversely, confining now our attention to the parabolic-elliptic-elliptic system  \eqref{problem} with $\rho=1$ (we just mention that other variants with nonlinear diffusion and/or chemo-sensitivities or logistic sources are available in the literature),  it is known that for $n\geq 2$ and for  $\xi\gamma-\chi\alpha>0$ (repulsion dominates attraction)  it only admits globally bounded classical solutions, whilst in the bi-dimensional setting and for $\xi\gamma-\chi\alpha<0$  unbounded solutions can be detected (see \cite{EspejoSuzukiAttractionRepulsion,TaoWangAttractionRepulsion,YuGuoZhengAttractionRepulsion}).  
Additionally, in \cite{CriticalMassAttrRepulGuoEtAl} the authors establish that for $n=2$ and $\chi\alpha-\xi\gamma>0$,  the value $\frac{4\pi}{\chi\alpha-\xi\gamma}$ is the critical mass deciding whether all solutions are global or, on the contrary, certain ones may blow-up; if this last scenario occurs, lower bounds of the blow-up time are estimated in \cite{Viglialoro2019RepulsionAttraction}.  
Now, by virtue of all of the considerations, we rigorously formulate  our main result. 
\begin{theorem}\label{Maintheorem}
Let $\Omega$ be a bounded and smooth domain of $\R^n$, $n\geq 2$. Then, for any $\alpha,\beta,\gamma,\delta,\chi,\xi>0$, $0<\rho<1 $, and any nonnegative and nontrivial initial data $0 \leq u_0(x)\in C^0(\bar{\Omega})$, problem \eqref{problem} admits a unique solution $(u,v,w)$ of nonnegative and bounded functions in the class
\begin{equation*}
C^0([0,\infty);C^0(\Omega))\cap C^{2,1}(\bar{\Omega}\times (0,\infty))\times C^{2,0}(\bar{\Omega}\times (0,\infty))\times C^{2,0}(\bar{\Omega}\times (0,\infty)).
\end{equation*}
\end{theorem}
\section{From local to globally bounded solutions}
One of the first steps in dealing with solutions of \eqref{problem} is showing that they do exist, at least locally.
\begin{lemma}\label{LocalExistenceLemma}  
Let $\Omega$ be a bounded and smooth domain of $\R^n$, $n\geq 2$. Then, for any $\alpha,\beta,\gamma,\delta\chi,\xi,\rho>0$, and any nonnegative and nontrivial initial data $0 \leq u_0(x)\in C^0(\bar{\Omega})$, problem \eqref{problem} admits a unique solution $(u,v,w)$ of nonnegative functions in the class
\begin{equation*}
C^0([0,T_{max});C^0(\Omega))\cap C^{2,1}(\bar{\Omega}\times (0,T_{max}))\times C^{2,0}(\bar{\Omega}\times (0,T_{max}))\times C^{2,0}(\bar{\Omega}\times (0,T_{max})).
\end{equation*}
Here $T_{max}\in(0,\infty]$, denoting the maximal existence time, is such that (dichotomy criterion) or $T_{max}=\infty$  (global-in-time classical solution) or if $T_{max}<\infty$ (local-in-time classical solution) then necessarily 
\begin{equation}\label{extensibility_criterion_Eq} 
\limsup_{t\nearrow T_{max}}\lVert u (\cdot,t)\rVert_{L^\infty(\Omega)}=\infty.
\end{equation}
Moreover, 
\begin{equation}\label{Bound_of_u} 
\int_\Omega u (\cdot ,t)  = m:=\int_\Omega u_0>0\quad \textrm{for all}\quad  t\in (0,T_{max}).
\end{equation}
\begin{proof}
The first statements can be shown by straightforward
 adaptations of well-established methods involving an appropriate fixed point framework and standard parabolic and elliptic regularity theory (see, for instance, 
\cite[Lemma 2.1]{FujieWinklerYokotaSignalSensi}), as well as related comparison principles. On the other hand, relation \eqref{Bound_of_u}  directly comes by integrating over $\Omega$  the equation for $u$ in \eqref{problem}. 
\end{proof}
\end{lemma}
Once the solvability (at least in the local sense) for problem \eqref{problem} is ensured, the bridge establishing the globability and boundedness  is achieved throughout some precise $L^p$-bound for these solutions. To be precise we have this
\begin{lemma}\label{FromLocalToGLobalBoundedLemma}
Under the assumptions of Lemma \ref{LocalExistenceLemma}, let  $(u,v,w)$  be the classical solution to problem \eqref{problem}. If for some $\frac{n}{2}<p<n$ the $u$-component belongs  to $L^\infty((0,T_{max});L^p(\Omega))$, then  $T_{max}=\infty$ and $u$ is uniformly bounded in $\Omega \times (0,\infty)$. 
\begin{proof}
Well-known elliptic regularity theory in conjunction with  Sobolev embedding theorems infer, through the second equation of \eqref{problem} and  $u\in L^\infty((0,T_{max});L^p(\Omega))$, that $v\in L^\infty((0,T_{max});W^{2,p}(\Omega))$, $\nabla v\in  L^\infty((0,T_{max});W^{1,p}(\Omega))$, and finally $ v\in  L^\infty((0,T_{max});C^{[2-(n/p)]}(\bar{\Omega}))$ and $\nabla v\in  L^\infty((0,T_{max});L^{q}(\Omega))$ for all $n<q<p^*:=\frac{np}{n-p}$. In particular, since the same reasoning is valid for $w$, by posing $\tilde{v}=\chi  u-\xi  w$ we have that for some positive constant $C_q$ 
\begin{equation}\label{Bound_v_1-q}
\lVert  \tilde{v} (\cdot, t)\lVert_{L^{q}(\Omega)}+\lVert  \nabla \tilde{v} (\cdot, t)\lVert_{L^{q}(\Omega)}  \leq C_q\quad \textrm{for all}\quad t\in(0,T_{max}).
\end{equation}
Additionally, for any $(x,t)\in \Omega \times (0,T_{max})$, the first equation of  \eqref{problem} reads $u_t=\Delta u-\nabla \cdot (u \nabla \tilde{v})$ so that for $t_0:=\max\{0,t-1\}$ the representation formula yields
\begin{equation*}
\begin{split}
u (\cdot,t) &\leq e^{(t-t_0)\Delta}u(\cdot,t_0)-\int_{t_0}^t e^{(t-s)\Delta}\nabla \cdot (u (\cdot,s) \nabla \tilde{v} (\cdot,s))ds=: u_{1}(\cdot,t)+u_{2}(\cdot,t). 
\end{split}
\end{equation*}
In these circumstances, the rest of the proof follows that done in \cite[Lemma 4.1]{ViglialoroWoolleyApplicableanalysis}; precisely,  in order to control the $L^\infty(\Omega)$-norm of $u$ on $(0,T_{max})$, first one can control a suitable  norm of the cross diffusion term $u \nabla \tilde{v}$ by replacing relation (24) therein with bound \eqref{Bound_v_1-q}, then applications of known smoothing estimates for the Neumann heat semigroup entail such uniform bound. Finally, the conclusion is achieved by relying on the dichotomy criterion  \eqref{extensibility_criterion_Eq}. 
\end{proof}
\end{lemma}
\begin{remark}
From the above lemma, the open question given in \cite[Remark 1]{Viglialoro2019RepulsionAttraction} has a response: indeed, in the context of \cite[Theorem 3.1]{Viglialoro2019RepulsionAttraction}, if the $u$-component of the solution $(u,v,w)$ to the bi-dimensional version of problem \eqref{problem} with $\rho=1$ becomes unbounded at some finite time $t^*$ (in the sense of the $L^\infty(\Omega)$-norm) it also blows-up in the $L^p(\Omega)$-norm for any $p>1$, since otherwise from Lemma \ref{FromLocalToGLobalBoundedLemma} with $n=2$ it would be globally bounded. In particular, that theorem continues valid also without the extra assumption that $E(t):=\int_\Omega u^p \nearrow \infty$ as $t\nearrow t^*$, therein required. 
\end{remark}
%
%
%
%
\section{Some properties of classical solutions: proof of the main theorem}
\subsection{Deriving a proper absorptive differential inequality}
With the crucial implication of Lemma \ref{FromLocalToGLobalBoundedLemma} in our hands, in this section we aim at bounding on $(0,T_{max})$ the functional $E(t):=\int_\Omega u^p$, for $p>1$, by means of a time independent constant. This will be obtained by deriving a proper absorptive differential inequality for $E$, exactly with the aid of this sequel of lemmas. 

%
%

\begin{lemma}\label{EllipticEhrlingSystemLemma}
Under the assumptions of Lemma \ref{LocalExistenceLemma}, let $(u,v,w)$ be the classical solution to problem \eqref{problem} and $m:=\int_\Omega u_0$. Then for any $p>1$ and $\sigma>0$ there exists $\tilde{c}=\tilde{c}(p,\sigma) >0$ such that the $w$-component 
satisfies, for any $\hat{c}>0,$ 
\begin{equation}\label{EhrlingTypeInequality}
\hat{c}\int_\Omega w^{p+1} \leq \sigma \int_\Omega u^{p+1}+\tilde{c}m^{p+1}\quad \textrm{for all } t \in (0,T_{max}) 
.
\end{equation}
Additionally, for $0<\theta=\frac{\frac{p}{2}-\frac{1}{2}}{\frac{p}{2}+\frac{1}{n}-\frac{1}{2}}<1$ and some constant $C_{GN}>0$, the $u$-component  fulfills 
\begin{equation}\label{InequalityTipoG-N} 
\int_\Omega u^p \leq \frac{4(p-1)}{p}\int_\Omega \lvert \nabla u^\frac{p}{2}\rvert^2 +c^* \quad \textrm{for all } t \in (0,T_{max}), 
\end{equation}
being $c^*=2m^pC_{GN}^2[(1-\theta)m^p(\frac{2(p-1)}{p \theta C_{GN}^2})^\frac{\theta}{\theta-1}+1].$
\begin{proof}
As to the first conclusion, for the sake of completeness and clarity we retrace in detail what presented in \cite[Lemma 2.2]{WinklerHowFar} and \cite[Lemma 2.2]{LankeitChemoPrevent}. A direct integration over $\Omega$ of the third equation in \eqref{problem} produces, for any $p>1$, and using \eqref{Bound_of_u} 
\begin{equation}\label{EqualityComingConservationMass}
\Big(\int_\Omega w \Big)^{p+1}  =\Big(\frac{\gamma}{\delta}\Big)^{p+1}\Big( \int_\Omega u \Big)^{p+1} =\Big(\frac{\gamma}{\delta}\Big)^{p+1}m^{p+1}\quad \textrm{ on } (0,T_{max}),
\end{equation}
whilst testing procedures and Young's inequality on the same equation yield
\begin{equation*}
\begin{split}
p\int_\Omega w^{p-1}\lvert \nabla w \rvert^2 &+\delta \int_\Omega w^{p+1} =\gamma \int_\Omega w^p u\leq \frac{4p}{(p+1)^2}\int_\Omega w^{p+1}  \\ & +\frac{\gamma^{p+1}}{4^p}(p+1)^{p-1}\int_\Omega u^{p+1}\quad \textrm{ for all } (0,T_{max}).
\end{split}
\end{equation*}
This, through the identity $ \lvert \nabla w^\frac{p+1}{2} \rvert^2=\frac{(p+1)^2}{4}w^{p-1} \lvert \nabla w \rvert^2$, reads for all $\eta \in (0,\frac{1}{2})$
\begin{equation}\label{InequalityLikeEhrling_1}
\eta \int_\Omega \lvert \nabla w^\frac{p+1}{2} \rvert^2 \leq \eta\int_\Omega w^{p+1}  +\eta \frac{\gamma^{p+1}}{4^{p+1}p}(p+1)^{p+1}\int_\Omega u^{p+1} \quad \textrm{ on } (0,T_{max}).
\end{equation}
On the other hand, for the same $\eta \in (0,\frac{1}{2})$, by virtue of the inclusions 
\[W^{1,2}(\Omega)\hookrightarrow \hookrightarrow L^2(\Omega)\hookrightarrow L^\frac{2}{p+1}(\Omega),\]
Ehrling's Lemma (see \cite[Lemma 1.1]{ShowalterBOOK}) yields a constant $c_E(\eta)>0$ such that 
\begin{equation*}
\|V\|_{L^2(\Omega)}^2\leq \eta \|V\|_{W^{1,2}(\Omega)}^2+c_E(\eta)\|V\|_{L^\frac{2}{p+1}(\Omega)}^2\quad \textrm{for all } V \in W^{1,2}(\Omega);
\end{equation*}
subsequently, posing in this last relation $V=w^\frac{p+1}{2}$, and making use of \eqref{EqualityComingConservationMass} and \eqref{InequalityLikeEhrling_1},  as well as of the conservation of mass property \eqref{Bound_of_u}, we obtain for all $t\in (0,T_{max})$
\begin{equation*}
(1-2\eta)\int_\Omega w^{p+1} \leq \eta \frac{(\gamma(p+1))^{p+1}}{4^{p+1}p}\int_\Omega u^{p+1} +\Big(\frac{\gamma}{\delta}\Big)^{p+1}c_E(\eta) m^{p+1}.
\end{equation*}
Finally, for any $\hat{c}>0$ we introduce the function $\sigma: (0,\frac{1}{2})\rightarrow (0,\infty)$ defined as $\sigma(\eta)=\frac{\eta}{1-2\eta}\frac{(\gamma(p+1))^{p+1}\hat{c}}{4^{p+1}p}$, and estimate \eqref{EhrlingTypeInequality} follows with the choice $$\tilde{c}=\tilde{c}(p,\sigma):=\big(\frac{\gamma}{\delta}\big)^{p+1}\Big(\hat{c}+\frac{2\sigma}{\frac{(\gamma(p+1))^{p+1}\hat{c}}{4^{p+1}p}}\Big)c_E\Big(\frac{\sigma}{2\sigma+\frac{(\gamma(p+1))^{p+1}\hat{c}}{4^{p+1}p}}\Big).$$ 

In turn, the proof of \eqref{InequalityTipoG-N} comes from an application of a general case of the Gagliardo--Nirenberg inequality: in particular, for any $p>1$, we can use \cite[(22) of Lemma 4]{ViglialoroApplMathOpt2019} with $f=u^\frac{p}{2}$, $\mathfrak{p}=\mathfrak{q}=2$ and $\mathfrak{r}=\frac{2}{p}$ so to explicitly have  
\begin{equation*}
\begin{split}
\int_\Omega u^p&= \|  u^\frac{p}{2} \|_{L^{2}(\Omega)}^2 \leq C_{GN}^2 ( \| \nabla u^\frac{p}{2} \|_{L^{2}(\Omega)}^{\theta} \| u^\frac{p}{2} \|_{L^{\frac{2}{p}}(\Omega)}^{1 - \theta}+  \| u^\frac{p}{2} \|_{L^{\frac{2}{p}}(\Omega)})^2 \\ &
= 2C_{GN}^2m^{p(1-\theta)}\Big(\int_\Omega |\nabla u^ \frac{p}{2}|^2\Big)^\theta+2C_{GN}^2m^{p}\quad \textrm{ on } (0,T_{max}); 
\end{split}
\end{equation*}
hence, we conclude invoking the Young inequality with exponents $\theta$ and $(1-\theta).$ (We remark that the proof of inequality \eqref{InequalityTipoG-N} does not rely on the fact that $(u,v,w)$ solves \eqref{problem}. Such an estimate, indeed, holds for any general function belonging to $W^{1,2}(\Omega)\cap L^\frac{2}{p}(\Omega)$; nevertheless,  to facilitate the reading we preferred to present it in this way.)
%
%
%
%
\end{proof}
\end{lemma}
\begin{lemma}\label{MainInequalityU^2Lemma}
Under the assumptions of Lemma \ref{LocalExistenceLemma}, but for $0<\rho<1$, let $(u,v,w)$ be the classical solution to problem \eqref{problem} and $m:=\int_\Omega u_0$. Then for any $p>1$  and $\tilde{c}=\tilde{c}(p,\frac{\gamma\xi(p-1)}{3})$ taken from Lemma \ref{EllipticEhrlingSystemLemma}, the $u$-component satisfies 
\begin{equation}\label{MainInequalityU^2}
\frac{d}{d t}\int_\Omega u^p \leq -\frac{4(p-1)}{p}\int_\Omega \lvert \nabla u^\frac{p}{2} \rvert^2+\bar{c} \quad \textrm{ for all} \quad  t \in (0,T_{max}),
\end{equation}
being $\bar{c}=c_1+\tilde{c}m^{p+1}$, with $c_1=\frac{\alpha\chi(p-1)(1-\rho)}{p+1}(\frac{(p+1)\gamma\xi}{(p+\rho)3\alpha \chi})^\frac{p+\rho}{\rho-1}|\Omega|$.
 \begin{proof}
For any $p>1$ by using problem \eqref{problem} and the divergence theorem (this in particular twice in both cross-diffusion terms), we have for all $t\in (0,T_{max})$
\begin{equation}\label{DerivativeE-FirstStep}
\begin{split}
\frac{d}{dt}\int_\Omega u^p&=p\int_\Omega u^{p-1} u_t  =p\int_\Omega u^{p-1} [\Delta u -\chi \nabla \cdot (u\nabla v)+\xi \nabla \cdot (u\nabla w) ]\\&
=-p(p-1)\int_\Omega u^{p-2} |\nabla u|^2-\chi\beta (p-1)\int_\Omega u^p v \\& \quad+ \alpha\chi (p-1)\int_\Omega u^{p+\rho} + \xi\delta (p-1)\int_\Omega u^pw-\gamma\xi (p-1)\int_\Omega u^{p+1}.
\end{split}
\end{equation}
On the other hand, if we neglect the nonpositive term $-\chi\beta (p-1)\int_\Omega u^p v $, use the Young inequality and  \eqref{EhrlingTypeInequality} with  $\tilde{c}$  as in our hypotheses (corresponding to the choice $\sigma=\frac{\gamma\xi(p-1)}{3}$),  $c_1$ as established in this statement  and $\hat{c}=\frac{\xi\delta(p-1)}{p+1}(\frac{(p+1)\gamma}{3p\delta})^{-p}$, we have that for all $t\in(0,T_{max})$ these relations are complied:
\begin{equation*}
\xi\delta (p-1)\int_\Omega u^pw \leq \frac{\gamma\xi(p-1)}{3}\int_\Omega u^{p+1}+\hat{c}\int_\Omega w^{p+1} \leq \frac{2\gamma\xi(p-1)}{3}\int_\Omega u^{p+1}+\tilde{c}m^{p+1},
\end{equation*}
and
\begin{equation*}
 \alpha\chi (p-1)\int_\Omega u^{p+\rho}\leq \frac{\gamma\xi(p-1)}{3}\int_\Omega u^{p+1}+c_1
 \quad \textrm{ on } (0,T_{max}).
\end{equation*}
Then, by virtue of the pointwise identity $u^{p-2}|\nabla u|^2=\frac{4}{p^2}|\nabla u^\frac{p}{2}|^2$ and the previous two inequalities, estimate
\eqref{DerivativeE-FirstStep} actually reads as claim \eqref{MainInequalityU^2}, with $\bar{c}=c_1+\tilde{c}m^{p+1}.$
\end{proof}
\end{lemma}
\begin{lemma}\label{MainfinalLemma}
Under the assumptions of Lemma \ref{LocalExistenceLemma},  but for $0<\rho<1$, let $(u,v,w)$ be the classical solution to problem \eqref{problem}. Then for all $p>1$ there exists $C>0$ such that 
\begin{equation*}
\int_\Omega u^p \leq C \quad \textrm{ for all} \quad  t \in (0,T_{max}).
\end{equation*}
 \begin{proof}
Collecting \eqref{InequalityTipoG-N}  and \eqref{MainInequalityU^2}, for $E(t)=\int_\Omega u^p$ and $c_{*}=c^{*}+\bar{c}$, we get the absorptive differential inequality $E'(t)\leq c_{*}- E(t)$ on  $(0,T_{max})$ which, complemented with the natural initial condition $E(0)=\int_\Omega u_0^p$, manifestly leads to $E(t)\leq \max\{E(0),c_{*}\}=:C$, for all $t\in(0,T_{max}).$
 \end{proof}
 \end{lemma}
\subsection{Proof of Theorem \ref{Maintheorem}}\label{MainTheoremProofSection} 
 For any $n\geq 2$, let $(u,v,w)$ be the classical solution  to \eqref{problem} provided by Lemma \ref{LocalExistenceLemma}. By choosing $\frac{n}{2}<p<n$, we have  from Lemma \ref{MainfinalLemma} that $u\in L^\infty((0,T_{max});L^p(\Omega))$, so that in turn an application of  Lemma \ref{FromLocalToGLobalBoundedLemma}
immediately concludes the proof. 
 \qed



\begin{thebibliography}{10}

\bibitem{EspejoSuzukiAttractionRepulsion}
E.~Espejo and T.~Suzuki.
\newblock Global existence and blow-up for a system describing the aggregation
  of microglia.
\newblock {\em Appl. Math. Lett.}, 35:29--34, 2014.

\bibitem{FujieWinklerYokotaSignalSensi}
K.~Fujie, M.~Winkler, and T.~Yokota.
\newblock Boundedness of solutions to parabolic-elliptic {K}eller-{S}egel
  systems with signal-dependent sensitivity.
\newblock {\em Math. Methods Appl. Sci.}, 38(6):1212--1224, 2015.

\bibitem{CriticalMassAttrRepulGuoEtAl}
Q.~Guo, Z.~Jiang, and S.~Zheng.
\newblock Critical mass for an attraction--repulsion chemotaxis system.
\newblock {\em Appl. Anal.}, 97(13):2349--2354, 2018.

\bibitem{K-S-1970}
E.~F. Keller and L.~A. Segel.
\newblock Initiation of slime mold aggregation viewed as an instability.
\newblock {\em J. Theor. Biol.}, 26(3):399--415, 1970.

\bibitem{Keller-1971-MC}
E.~F. Keller and L.~A. Segel.
\newblock {Model for chemotaxis.}
\newblock {\em J. Theor. Biol.}, 30(2):225--234, 1971.

\bibitem{LankeitChemoPrevent}
J.~Lankeit.
\newblock Chemotaxis can prevent thresholds on population density.
\newblock {\em Discrete Continuous Dyn. Syst. Ser. B.}, 20(5):1499--1527, 2015.

\bibitem{Luca2003Alzheimer}
M.~Luca, A.~Chavez-Ross, L.~Edelstein-Keshet, and A.~Mogilner.
\newblock {Chemotactic signaling, microglia, and Alzheimer's disease senile
  plaques: Is there a connection?}
\newblock {\em Bull. Math. Biol.}, 2003.

\bibitem{ShowalterBOOK}
R.~E. Showalter.
\newblock {\em Monotone Operators in Banach Space and Nonlinear Partial
  Differential Equations}.
\newblock American Mathematical Society, 1997.

\bibitem{TaoWangAttractionRepulsion}
Y.~Tao and Z.-A. Wang.
\newblock Competing effects of attraction vs. repulsion in chemotaxis.
\newblock {\em Math. Models Methods Appl. Sci.}, 23(01):1--36, 2013.

\bibitem{ViglialoroApplMathOpt2019}
G.~Viglialoro.
\newblock Global in time and bounded solutions to a parabolic-elliptic
  chemotaxis system with nonlinear diffusion and signal-dependent sensitivity.
\newblock {\em Appl. Math. Opt.}, 2019.
  \href{https://link.springer.com/article/10.1007/s00245-019-09575-0}{doi:
  10.1007/s00245-019-09575-0}.

\bibitem{Viglialoro2019RepulsionAttraction}
G.~Viglialoro.
\newblock Explicit lower bound of blow-up time for an attraction-repulsion
  chemotaxis system.
\newblock {\em J. Math. Anal. Appl.}, 2019.
  \href{https://www.sciencedirect.com/science/article/pii/S0022247X19305414}{doi:
  10.1016/j.jmaa.2019.06.067}.

\bibitem{ViglialoroWoolleyApplicableanalysis}
G.~Viglialoro and T.~Woolley.
\newblock {Solvability of a Keller--Segel system with signal-dependent
  sensitivity and essentially sublinear production}.
\newblock {\em Appl. Anal.}, 2019.
  \href{https://www.tandfonline.com/doi/full/10.1080/00036811.2019.1569227}{doi:
  10.1080/00036811.2019.1569227}.

\bibitem{WinklerHowFar}
M.~Winkler.
\newblock How far can chemotactic cross-diffusion enforce exceeding carrying
  capacities?
\newblock {\em J. Nonlinear. Sci.}, 24(5):809--855, 2014.

\bibitem{YuGuoZhengAttractionRepulsion}
H.~Yu, Q.~Guo, and S.~Zheng.
\newblock Finite time blow-up of nonradial solutions in an
  attraction--repulsion chemotaxis system.
\newblock {\em Nonlinear Anal. Real. World Appl.}, 34:335--342, 2017.

\end{thebibliography}
\end{document}